\documentclass{article}
\usepackage[centertags]{amsmath}
\usepackage{amsfonts}
\usepackage{amssymb}
\usepackage{amsthm}
\usepackage{newlfont}
\newcommand{\al}{\alpha}
\newcommand{\Ri}{\mathrm{Ric}}
\newcommand{\co}{\mathrm{const}}
\newcommand{\e}{\epsilon}
\newcommand{\de}{\delta}
\newcommand{\Om}{\Omega}
\newcommand{\be}{\beta}
\newcommand{\ka}{\kappa}
\newcommand{\la}{\lambda}
\newcommand{\mi}{\mathrm{min}}
\numberwithin{equation}{section}
\begin{document}
\title{Ricci flow with surgery on three-manifolds}\author{Grisha
Perelman\thanks{St.Petersburg branch of Steklov Mathematical
Institute, Fontanka 27, St.Petersburg 191011, Russia. Email:
perelman@pdmi.ras.ru or perelman@math.sunysb.edu }} \maketitle
\par This is a technical paper, which is a continuation of [I].
Here we verify most of the assertions, made in [I, \S 13]; the
exceptions are  (1) the statement that a 3-manifold which
collapses with local lower bound for sectional curvature is a
graph manifold - this is deferred to a separate paper, as the
proof has nothing to do with the Ricci flow, and (2) the claim
about the lower bound for the volumes of the maximal horns and the
smoothness of the solution from some time on, which turned out to
be unjustified, and, on the other hand, irrelevant for the other
conclusions.
\par The Ricci flow with surgery was considered by Hamilton [H 5,\S
4,5]; unfortunately, his argument, as written, contains an
unjustified statement ($R_{MAX}=\Gamma,$ on page 62, lines 7-10
from the bottom), which I was unable to fix. Our approach is
somewhat different, and is aimed at eventually constructing a
canonical Ricci flow, defined on a largest possible subset of
space-time, - a goal, that has not been achieved yet in the
present work. For this reason, we consider two scale bounds: the
cutoff radius $h,$ which is the radius of the necks, where the
surgeries are performed, and the much larger radius $r,$ such that
the solution on the scales less than $r$ has standard geometry.
The point is to make $h$ arbitrarily small while keeping $r$
bounded away from zero.
\section*{Notation and terminology}
\par $\ \ \ \ \ \ B(x,t,r)$ denotes the open metric ball of radius $r,$ with respect to the
metric at time $t,$ centered at $x.$
\par $P(x,t,r,\triangle t)$ denotes a parabolic neighborhood, that
is the set of all points $(x',t')$ with $x'\in B(x,t,r)$ and
$t'\in [t,t+\triangle t]$ or $t'\in [t+\triangle t,t],$ depending
on the sign of $\triangle t.$
\par A ball $B(x,t,\e^{-1}r)$ is called an $\e$-neck, if, after scaling
the metric with factor $r^{-2},$ it is $\e$-close to the standard
neck $\mathbb{S}^2\times\mathbb{I},$ with the product metric,
where $\mathbb{S}^2$ has constant scalar curvature one, and
$\mathbb{I}$ has length $2\e^{-1};$ here $\e$-close refers to
$C^N$ topology, with $N>\e^{-1}.$
\par A parabolic neighborhood $P(x,t,\e^{-1}r,r^2)$ is called a
strong $\e$-neck, if, after scaling with factor $r^{-2},$ it is
$\e$-close to the evolving standard neck, which at each time
$t'\in [-1,0]$ has length $2\e^{-1}$ and scalar curvature
$(1-t')^{-1}.$
\par A metric on $\mathbb{S}^2\times\mathbb{I},$ such that each point
is contained in some $\e$-neck, is called an $\e$-tube, or an
$\e$-horn, or a double $\e$-horn, if the scalar curvature stays
bounded on both ends, stays bounded on one end and tends to
infinity on the other, and tends to infinity on both ends,
respectively.
\par A metric on $\mathbb{B}^3$ or
$\mathbb{RP}^3\setminus\bar{\mathbb{B}}^3,$ such that each point
outside some compact subset is contained in an $\e$-neck, is
called an $\e$-cap or a capped $\e$-horn, if the scalar curvature
stays bounded or tends to infinity on the end, respectively.
\par We denote by $\e$ a fixed small positive constant.
In contrast, $\de$ denotes a positive quantity, which is supposed
to be as small as needed in each particular argument.
\section{Ancient solutions with bounded entropy}
\par {\bf 1.1} In this section we review some of the results, proved or
quoted in [I,\S 11], correcting a few inaccuracies. We consider
smooth solutions $g_{ij}(t)$ to the Ricci flow on oriented
3-manifold $M$, defined for $-\infty<t\le 0$, such that for each
$t$ the metric $g_{ij}(t)$ is a complete non-flat metric of
bounded nonnegative sectional curvature, $\kappa$-noncollapsed on
all scales for some fixed $\kappa>0;$ such solutions will be
called ancient $\kappa$-solutions for short. By Theorem I.11.7,
the set of all such solutions with fixed $\kappa$ is compact
modulo scaling, that is from any sequence of such solutions
$(M^{\al},g_{ij}^{\al}(t))$ and points $(x^{\al},0)$ with
$R(x^{\al},0)=1$, we can extract a smoothly (pointed) convergent
subsequence, and the limit $(M,g_{ij}(t))$ belongs to the same
class of solutions. (The assumption in I.11.7. that $M^{\al}$ be
noncompact was clearly redundant, as it was not used in the proof.
Note also that $M$ need not have the same topology as $M^{\al}.)$
Moreover, according to Proposition I.11.2, the scalings of any
ancient $\kappa$-solution $g_{ij}(t)$ with factors $(-t)^{-1}$
about appropriate points converge along a subsequence of $t\to
-\infty$ to a non-flat gradient shrinking soliton, which will be
called an asymptotic soliton of the ancient solution. If the
sectional curvature of this asymptotic soliton is not strictly
positive, then by Hamilton's strong maximum principle it admits
local metric splitting, and it is easy to see that in this case
the soliton is either the round infinite cylinder, or its
$\mathbb{Z}_2$ quotient, containing one-sided projective plane. If
the curvature is strictly positive and the soliton is compact,
then it has to be a metric quotient of the round 3-sphere, by [H
1]. The noncompact case is ruled out below.
\par {\bf 1.2 Lemma.} {\it There is no (complete oriented
3-dimensional) noncompact $\kappa$-noncollapsed gradient shrinking
soliton with bounded positive sectional curvature.}
\par {\it Proof.} A gradient shrinking soliton $g_{ij}(t),
-\infty<t<0,$ satisfies the equation \begin{equation}
\nabla_i\nabla_j f+R_{ij}+\frac{1}{2t}g_{ij}=0  \end{equation}
Differentiating and switching the order of differentiation, we get
\begin{equation} \nabla_i R=2R_{ij}\nabla_j f \end{equation}
Fix some $t<0,$ say $t=-1,$ and consider a long shortest geodesic
$\gamma(s), 0\le s\le \bar{s};$ let $x=\gamma(0),
\bar{x}=\gamma(\bar{s}), X(s)=\dot{\gamma}(s).$ Since the
curvature is bounded and positive, it is clear from the second
variation formula that $ \int_0^{\bar{s}}{\Ri(X,X)ds}\le \co.$
Therefore, $\int_0^{\bar{s}}{|\Ri(X,\cdot)|^2ds}\le \co,$ and
$\int_0^{\bar{s}}{|\Ri(X,Y)|ds}\le\co (\sqrt{\bar{s}}+1)$ for any
unit vector field $Y$ along $\gamma,$ orthogonal to $X.$ Thus by
integrating (1.1) we get $X\cdot f(\gamma(\bar{s}))\ge
\frac{\bar{s}}{2}+\co, |Y\cdot f(\gamma(\bar{s}))|\le \co
(\sqrt{\bar{s}}+1).$ We conclude that at large distances from
$x_0$ the function $f$ has no critical points, and its gradient
makes small angle with the gradient of the distance function from
$x_0.$ \par Now from (1.2) we see that $R$ is increasing along the
gradient curves of $f,$ in particular, $\bar{R}=\mathrm{lim \ sup
\  }R>0.$ If we take a limit of our soliton about points
$(x^{\alpha}, -1)$ where $R(x^{\alpha})\to\bar{R},$ then we get an
ancient $\kappa$-solution, which splits off a line, and it follows
from I.11.3, that this solution is the shrinking round infinite
cylinder with scalar curvature $\bar{R}$ at time $t=-1.$  Now
comparing the evolution equations for the scalar curvature on a
round cylinder and for the asymptotic scalar curvature on a
shrinking soliton we conclude that $\bar{R}=1.$ Hence, $R(x)<1$
when the distance from $x$ to $x_0$ is large enough, and $R(x)\to
1$ when this distance tends to infinity.
\par Now let us check that the level surfaces of $f,$
sufficiently distant from $x_0,$ are convex. Indeed, if $Y$ is a
unit tangent vector to such a surface, then $\nabla_Y\nabla_Y f=
\frac{1}{2}-\Ri(Y,Y)\ge \frac{1}{2}-\frac{R}{2}>0.$ Therefore, the
area of the level surfaces grows as $f$ increases, and is
converging to the area of the round sphere of scalar curvature
one. On the other hand, the intrinsic scalar curvature of a level
surface turns out to be less than one. Indeed,  denoting by $X$
the unit normal vector, this intrinsic curvature can be computed
as $$ R-2\Ri(X,X)+2\frac{\mathrm{det}(\mathrm{Hess} f)}{|\nabla
f|^2}\le R-2\Ri(X,X)+\frac{(1-R+\Ri(X,X))^2}{2|\nabla f|^2}<1 $$
when $R$ is close to one and $|\nabla f|$ is large. Thus we get a
contradiction to the Gauss-Bonnet formula.
\par {\bf 1.3} Now, having listed all the asymptotic solitons, we
can classify the ancient $\kappa$-solutions. If such a solution
has a compact asymptotic soliton, then it is itself a metric
quotient of the round 3-sphere, because the positive curvature
pinching can only improve in time [H 1]. If the asymptotic soliton
contains the one-sided projective plane, then the solution  has a
$\mathbb{Z}_2$ cover, whose asymptotic soliton is the round
infinite cylinder. Finally, if the asymptotic soliton is the
cylinder,then the solution can be either noncompact (the round
cylinder itself, or the Bryant soliton, for instance), or compact.
The latter possibility, which was overlooked in the first
paragraph of [I.11.7], is illustrated by the example below, which
also gives the negative answer to the question in the very end of
[I.5.1].
\par {\bf 1.4 Example.} Consider a solution to the Ricci flow,
starting from a metric on $\mathbb{S}^3$ that looks like a long
round cylinder $\mathbb{S}^2\times\mathbb{I}$ (say, with radius
one and length $L>>1$), with two spherical caps, smoothly attached
to its boundary components. By [H 1] we know that the flow shrinks
such a metric to a point in time, comparable to one (because both
the lower bound for scalar curvature and the upper bound for
sectional curvature are comparable to one) , and after
normalization, the flow converges to the round 3-sphere. Scale the
initial metric and choose the time parameter in such a way that
the flow starts at time $t_0=t_0(L)<0,$ goes singular at $t=0,$
and at $t=-1$ has the ratio of the maximal sectional curvature to
the minimal one equal to $1+\e.$ The argument in [I.7.3] shows
that our solutions are $\kappa$-noncollapsed for some $\kappa>0$
independent of $L.$ We also claim that $t_0(L)\to -\infty$ as
$L\to\infty.$ Indeed, the Harnack inequality of Hamilton [H 3]
implies that $\ \ R_t\ge\frac{R}{t_0-t},\ \ \ \ $ hence $\ \
R\le\frac{2(-1-t_0)}{t-t_0}$ for $t\le -1,\ \ \ \ \ \ \ \ $ and
then the distance change estimate $\ \ \
\frac{d}{dt}\mathrm{dist}_t(x,y)\ge
-\co\sqrt{R_{\mathrm{max}}(t)}\ \ \ \ $ from  [H 2,\S 17] implies
that the diameter of $g_{ij}(t_0)$ does not exceed $ -\co\cdot
t_0$, which is less than $L\sqrt{-t_0}$ unless $t_0$ is large
enough. Thus, a subsequence of our solutions with $L\to\infty$
converges to an ancient $\kappa$-solution on $\mathbb{S}^3,$ whose
asymptotic soliton can not be anything but the cylinder.
\par {\bf 1.5} The important conclusion from the classification above and the
proof of Proposition I.11.2 is that there exists $\kappa_0>0,$
such that every ancient $\kappa$-solution is either
$\kappa_0$-solution, or a metric quotient of the round sphere.
Therefore, the compactness theorem I.11.7 implies the existence of
a universal constant $\eta,$ such that at each point of every
ancient $\kappa$-solution we have estimates \begin{equation}
|\nabla R|<\eta R^{\frac{3}{2}}, |R_t|<\eta R^2 \end{equation}
Moreover, for every sufficiently small $\e>0$ one can find
$C_{1,2}=C_{1,2}(\e),$ such that  for each point $(x,t)$ in every
ancient $\kappa$-solution there is a radius $r,
0<r<C_1R(x,t)^{-\frac{1}{2}},$ and a neighborhood $B,
B(x,t,r)\subset B\subset B(x,t,2r),$ which falls into one of the
four categories:
\par (a) $B$ is a strong $\e$-neck (more precisely, the slice of a strong $\e$-neck at its
maximal time), or
\par (b) $B$ is an $\e$-cap, or
\par (c) $B$ is a closed manifold, diffeomorphic to $\mathbb{S}^3$
or $\mathbb{RP}^3,$ or
\par (d) $B$ is a closed manifold of constant positive sectional
curvature;
\par furthermore, the scalar curvature in
$B$ at time $t$ is between $C_2^{-1}R(x,t)$ and  $C_2R(x,t),$ its
volume in cases (a),(b),(c) is greater than
$C_2^{-1}R(x,t)^{-\frac{3}{2}},$ and in case (c) the sectional
curvature in $B$ at time $t$ is greater than $C_2^{-1}R(x,t).$
\section{The standard solution}
\par Consider a rotationally symmetric metric on $\mathbb{R}^3$
with nonnegative sectional curvature, which splits at infinity as
the metric product of a ray and the round 2-sphere of scalar
curvature one. At this point we make some choice for the metric on
the cap, and will refer to it as the standard cap; unfortunately,
the most obvious choice, the round hemisphere, does not fit,
because the metric on $\mathbb{R}^3$ would not be smooth enough,
however we can make our choice as close to it as we like. Take
such a metric on $\mathbb{R}^3$ as the initial data for a solution
$g_{ij}(t)$ to the Ricci flow on some time interval $[0,T),$ which
has bounded curvature for each $t\in [0,T).$
\par {\bf Claim 1.} {\it The solution is rotationally symmetric
for all $t.$}
\par Indeed, if $u^i$ is a vector field evolving by
$u^i_t=\triangle u^i+R^i_ju^j,$ then $v_{ij}=\nabla_i u_j$ evolves
by $(v_{ij})_t=\triangle
v_{ij}+2R_{ikjl}v_{kl}-R_{ik}v_{kj}-R_{kj}v_{ik}.$ Therefore, if
$u^i$ was a Killing field at time zero, it would stay Killing by
the maximum principle. It is also clear that the center of the
cap, that is the unique maximum point for the Busemann function,
and the unique point, where all the Killing fields vanish, retains
these properties, and the gradient of the distance function from
this point stays orthogonal to all the Killing fields. Thus, the
rotational symmetry is preserved.
\par {\bf Claim 2.} {\it The solution converges at infinity to the
standard solution on the round infinite cylinder of scalar
curvature one. In particular, $T\le 1.$}
\par {\bf Claim 3.} {\it The solution is unique.}
\par Indeed, using Claim 1, we can reduce the linearized Ricci
flow equation to the system of two equations on
$(-\infty,+\infty)$ of the following type $$f_t=f''+a_1 f'+b_1
g'+c_1 f+d_1 g,\ \ g_t=a_2f'+b_2 g'+c_2 f+d_2 g,$$ where the
coefficients and their derivatives are bounded, and the unknowns
$f,g$ and their derivatives tend to zero at infinity by Claim 2.
So we get uniqueness by looking at the integrals
$\int_{-A}^A{(f^2+g^2)}$ as $A\to \infty.$
\par {\bf Claim 4.} {\it The solution can be extended to the time interval
$[0,1).$}
\par Indeed, we can obtain our solution as a limit of the
solutions on $\mathbb{S}^3,$ starting from the round cylinder
$\mathbb{S}^2\times\mathbb{I}$ of length $L$ and scalar curvature
one, with two caps attached; the limit is taken about the center
$p$ of one of the caps, $L\to\infty.$ Assume that our solution
goes singular at some time $T<1.$ Take $T_1<T$ very close to $T, \
\ T-T_1<<1-T.$ By Claim 2, given $\de>0,$ we can find
$\bar{L},\bar{D}<\infty,$ depending on $\de$ and $T_1,$ such that
for any point $x$ at distance $\bar{D}$ from $p$ at time zero, in
the solution with $L\ge\bar{L},$ the ball $B(x,T_1,1)$ is
$\de$-close to the corresponding ball in the round cylinder of
scalar curvature $(1-T_1)^{-1}.$ We can also find $r=r(\de ,T),$
independent of $T_1,$ such that the ball $B(x,T_1,r)$ is
$\de$-close to the corresponding euclidean ball. Now we can apply
Theorem I.10.1 and get a uniform estimate on the curvature at $x$
as $t\to T$, provided that $T-T_1<\e^2 r(\de,T)^2.$ Therefore, the
$t\to T$ limit  of our limit solution on the capped infinite
cylinder will be smooth near $x.$ Thus, this limit will be a
positively curved space with a conical point. However, this leads
to a contradiction via a blow-up argument; see the end of the
proof of the Claim 2 in I.12.1.
\par The solution constructed above will be called the standard
solution.
\par {\bf Claim 5.} {\it The standard solution satisfies the
conclusions of 1.5 , for an appropriate choice of  $\e, \ \eta ,
C_1(\e),C_2(\e),$ except that the $\e$-neck neighborhood need not
be strong; more precisely, we claim that if $(x,t)$ has neither an
$\e$-cap neighborhood as in 1.5(b), nor a strong $\e$-neck
neighborhood as in 1.5(a), then $x$ is not in $B(p,0,\e^{-1}),\ \
t<3/4,$ and there is an $\e$-neck $B(x,t,\e^{-1}r),$ such that the
solution in $P(x,t,\e^{-1}r,-t)$ is, after scaling with factor
$r^{-2}, \ \ \e$-close to the appropriate piece of the evolving
round infinite cylinder.
\par Moreover, we have an estimate $R_{\mathrm{min}}(t)\ge \co\cdot
(1-t)^{-1}.$ }
\par Indeed, the
statements follow from compactness  and Claim 2 on compact
subintervals of $[0,1),$ and from the same arguments as for
ancient solutions, when $t$ is close to one.
\section{The structure of solutions at the first singular time}
\par Consider a smooth solution $g_{ij}(t)$ to the Ricci flow on
$M\times [0,T),$ where $M$ is a closed oriented 3-manifold,
$T<\infty.$ Assume that curvature of $g_{ij}(t)$ does not stay
bounded as $t\to T.$ Recall that we have a pinching estimate
$Rm\ge -\phi(R)R$ for some function $\phi$ decreasing to zero at
infinity [H 4,\S 4], and that the solution is
$\kappa$-noncollapsed on the scales $\le r$ for some $\kappa>0,
r>0$ [I, \S 4].Then by Theorem I.12.1 and the conclusions of 1.5
we can find $r=r(\e)>0,$ such that each point $(x,t)$ with
$R(x,t)\ge r^{-2}$ satisfies the estimates (1.3) and has a
neighborhood, which is either an $\e$-neck, or an $\e$-cap, or a
closed positively curved manifold. In the latter case the solution
becomes extinct at time $T,$ so we don't need to consider it any
more. \par If this case does not occur, then let $\Om$ denote the
set of all points in $M,$ where curvature stays bounded as $t\to
T.$ The estimates (1.3) imply that $\Om$ is open and that
$R(x,t)\to \infty$ as $t\to T$ for each $x\in M\backslash\Om.$ If
$\Om$ is empty, then the solution becomes extinct at time $T$ and
it is entirely covered by $\e$-necks and caps shortly before that
time, so it is easy to see that $M$ is diffeomorphic to either
$\mathbb{S}^3$, or $\mathbb{RP}^3$, or
$\mathbb{S}^2\times\mathbb{S}^1$, or $\mathbb{RP}^3\ \sharp\
\mathbb{RP}^3.$
 \par Otherwise, if $\Om$ is not empty, we may
(using the local derivative estimates due to W.-X.Shi, see [H 2,\S
13]) consider a smooth metric $\bar{g}_{ij}$ on $\Om ,$ which is
the limit of $g_{ij}(t)$ as $t\to T.$ Let $\Om_\rho$ for some
$\rho<r$ denotes the set of points $x\in\Om ,$ where the scalar
curvature $\bar{R}(x)\le \rho^{-2}.$ We claim that $\Om_\rho$ is
compact. Indeed, if $\bar{R}(x)\le \rho^{-2},$ then we can
estimate the scalar curvature $R(x,t)$ on $[T-\eta^{-1}\rho^2,T)$
using (1.3), and for earlier times by compactness, so $x$ is
contained in $\Om$ with a ball of definite size, depending on
$\rho.$
\par Now take any $\e$-neck in $(\Om,\bar{g}_{ij})$ and consider a
point $x$ on one of its boundary components. If $x\in
\Om\backslash\Om_{\rho} ,$ then there is either an $\e$-cap or an
$\e$-neck, adjacent to the initial $\e$-neck. In the latter case
we can take a point on the boundary of the second $\e$-neck and
continue. This procedure can either terminate when we reach a
point in $\Om_{\rho}$ or an $\e$-cap, or go on indefinitely,
producing an $\e$-horn. The same procedure can be repeated for the
other boundary component of the initial $\e$-neck. Therefore,
taking into account that $\Om$ has no compact components, we
conclude that each $\e$-neck of $(\Om ,\bar{g}_{ij})$ is contained
in a subset of $\Om$ of one of the following types:
\par (a) An $\e$-tube with boundary components in $\Om_{\rho} ,$ or
\par (b) An $\e$-cap with boundary in $\Om_{\rho} ,$ or
\par (c) An $\e$-horn with boundary in $\Om_{\rho} ,$ or
\par (d) A capped $\e$-horn, or
\par (e) A double $\e$-horn.
\par Clearly, each $\e$-cap, disjoint from $\Om_{\rho},$ is also
contained in one of the subsets above. It is also clear that there
is a definite lower bound (depending on $\rho$) for the volume of
subsets of types (a),(b),(c), so there can be only finite number
of them. Thus we can conclude that there is only a finite number
of components of $\Om,$ containing points of $\Om_{\rho},$ and
every such component has a finite number of ends, each being an
$\e$-horn. On the other hand, every component of $\Om,$ containing
no points of $\Om_{\rho},$ is either a capped $\e$-horn, or a
double $\e$-horn.
\par Now, by looking at our solution for times $t$ just before
$T,$ it is easy to see that the topology of $M$ can be
reconstructed as follows: take the components $\Om_j, 1\le j\le i$
of $\Om$ which contain points of $\Om_{\rho},$ truncate their
$\e$-horns, and glue to the boundary components of truncated
$\Om_j$ a collection of tubes $\mathbb{S}^2\times\mathbb{I}$ and
caps $\mathbb{B}^3$ or $\mathbb{RP}^3\backslash\mathbb{B}^3.$
Thus, $M$ is diffeomorphic to a connected sum of $\bar{\Om}_j ,
1\le j\le i,$ with a finite number of
$\mathbb{S}^2\times\mathbb{S}^1$ (which correspond to gluing a
tube to two boundary components of the same $\Om_j$), and a finite
number of $\mathbb{RP}^3;$ here $\bar{\Om}_j$ denotes $\Om_j$ with
each $\e$-horn one point compactified.
\section{Ricci flow with cutoff}
\par {\bf 4.1} Suppose we are given a collection of smooth
solutions $g_{ij}(t)$ to the Ricci flow, defined on
$M_k\times[t_k^-,t_k^+),$ which go singular as $t\to t_k^+.$ Let
$(\Om_k, \bar{g}_{ij}^k)$ be the limits of the corresponding
solutions as $t\to t_k^+,$ as in the previous section. Suppose
also that for each $k$ we have $t_k^- =t_{k-1}^+,$ and
$(\Om_{k-1},\bar{g}_{ij}^{k-1})$ and $(M_k,g_{ij}^k(t_k^-))$
contain compact (possibly disconnected) three-dimensional
submanifolds  with smooth boundary, which are isometric. Then we
can identify these isometric submanifolds and talk about the
solution to the Ricci flow with surgery on the union of all
$[t_k^-,t_k^+).$
\par Fix a small number $\e>0$ which is admissible in sections
1,2. In this section we consider only solutions to the Ricci flow
with surgery, which satisfy the following a priori assumptions:
\par (pinching) There exists a function $\phi,$ decreasing to zero at
infinity, such that $Rm\ge -\phi(R)R,$
\par (canonical neighborhood) There exists $r>0,$ such that every
point where scalar curvature is at least $r^{-2}$ has a
neighborhood, satisfying the conclusions of 1.5. (In particular,
this means that if in case (a) the neighborhood in question is
$B(x_0,t_0,\e^{-1}r_0),$ then the solution is required to be
defined in the whole $P(x_0,t_0,\e^{-1}r_0,-r_0^2);$ however, this
does not rule out a surgery in the time interval
$(t_0-r_0^2,t_0),$ that occurs sufficiently far from $x_0$.)
\par Recall that from the pinching estimate of Ivey and Hamilton, and Theorem
I.12.1, we know that the a priori assumptions above hold for a
smooth solution on any finite time interval. For Ricci flow with
surgery they will be justified in the next section.
\par {\bf 4.2 Claim 1.} {\it Suppose we have a solution to the Ricci
flow with surgery, satisfying the canonical neighborhood
assumption, and let $Q=R(x_0,t_0)+r^{-2}.$ Then we have  estimate
$R(x,t)\le 8Q$ for those $(x,t)\in
P(x_0,t_0,\frac{1}{2}\eta^{-1}Q^{-\frac{1}{2}},-\frac{1}{8}\eta^{-1}Q^{-1}),$
for which the solution is defined. }
\par Indeed, this follows from estimates (1.3).
\par {\bf Claim 2.} {\it For any $A<\infty$ one can find
$Q=Q(A)<\infty$ and $\xi=\xi(A)>0$ with the following property.
Suppose we have a solution to the Ricci flow with surgery,
satisfying the pinching and the canonical neighborhood
assumptions. Let $\gamma$ be a shortest geodesic in $g_{ij}(t_0)$
with endpoints $x_0$ and $x,$ such that $R(y,t_0)>r^{-2}$ for each
$y\in\gamma ,$ and $Q_0=R(x_0,t_0)$ is so large that
$\phi(Q_0)<\xi.$ Finally, let $z\in \gamma$ be any point
satisfying $R(z,t_0)>10C_2R(x_0,t_0).$ Then
$\mathrm{dist}_{t_0}(x_0,z)\ge AQ_0^{-\frac{1}{2}}$ whenever
$R(x,t_0)>QQ_0.$}
\par The proof is exactly the same as for Claim 2 in Theorem
I.12.1; in the very end of it, when we get a piece of a non-flat
metric cone as a blow-up limit, we get a contradiction to the
canonical neighborhood assumption, because the canonical
neighborhoods of types other than (a) are not close to a piece of
metric cone, and type (a) is ruled out by the strong maximum
principle, since the $\e$-neck in question is strong.
\par {\bf 4.3} Suppose we have a solution to the Ricci flow with
surgery, satisfying our a priori assumptions, defined on $[0,T),$
and going singular at time $T.$ Choose a small $\de>0$ and let
$\rho=\de r.$ As in the previous section, consider the limit
$(\Om,\bar{g}_{ij})$ of our solution as $t\to T,$ and the
corresponding compact set $\Om_{\rho}.$
\par {\bf Lemma.} {\it There exists a radius $h, 0<h<\de\rho,$
depending only on $\de ,\rho$ and the pinching function $\phi,$
such that for each point $x$ with
$h(x)=\bar{R}^{-\frac{1}{2}}(x)\le h$ in an $\e$-horn of
$(\Om,\bar{g}_{ij})$ with boundary in $\Om_{\rho},$ the
neighborhood $P(x,T,\de^{-1}h(x),-h^2(x))$ is a strong
$\de$-neck.}
\par {\it Proof.} An argument by contradiction. Assuming the
contrary, take a sequence of solutions with limit metrics
$(\Om^{\al},\bar{g}_{ij}^{\al})$ and points $x^{\al}$ with
$h(x^{\al})\to 0.$ Since $x^{\al}$ lies deeply inside an
$\e$-horn, its canonical neighborhood is a strong $\e$-neck. Now
Claim 2 gives  the curvature estimate that allows us to take a
limit of appropriate scalings of the metrics $g_{ij}^{\al}$ on
$[T-h^2(x^{\al}),T]$ about $x^{\al},$ for a subsequence of
$\al\to\infty.$ By shifting the time parameter we may assume that
the limit is defined on $[-1,0].$ Clearly, for each time in this
interval, the limit is a complete manifold with nonnegative
sectional curvature; moreover, since $x^{\al}$ was contained in an
$\e$-horn with boundary in $\Om^{\al}_{\rho},$ and
$h(x^{\al})/\rho\to 0,$ this manifold has two ends. Thus, by
Toponogov, it admits a metric splitting
$\mathbb{S}^2\times\mathbb{R}.$ This implies that the canonical
neighborhood of the point $(x^{\al},T-h^2(x^{\al}))$ is also of
type (a), that is a strong $\e$-neck, and we can repeat the
procedure to get the limit, defined on $[-2,0],$ and so on. This
argument works for the limit in any finite time interval $[-A,0],$
because $h(x^{\al})/\rho\to 0.$ Therefore, we can construct a
limit on $[-\infty,0];$ hence it is the round cylinder, and we get
a contradiction.
\par {\bf 4.4} Now we can specialize our surgery and define the
Ricci flow with $\de$-cutoff. Fix $\de>0,$ compute $\rho=\de r$
and determine $h$ from the lemma above. Given a smooth metric
$g_{ij}$ on a closed manifold, run the Ricci flow until it goes
singular at some time $t^+;$ form the limit $(\Om,\bar{g}_{ij}).$
If $\Om_{\rho}$ is empty, the procedure stops here, and we say
that the solution became extinct. Otherwise we remove the
components of $\Om$ which contain no points of $\Om_{\rho},$ and
in every $\e$-horn of each of the remaining components we find a
$\de$-neck of radius $h,$ cut it along the middle two-sphere,
remove the horn-shaped end, and glue in an almost standard cap in
such a way that the curvature pinching is preserved and a metric
ball of radius $(\de ')^{-1}h$ centered near the center of the cap
is, after scaling with factor $h^{-2},$ $\ \ \de '$-close to the
corresponding ball in the standard capped infinite cylinder,
considered in section 2. (Here $\de '$ is a function of $\de $
alone, which tends to zero with $\de.$)
\par The possibility of capping a $\de$-neck preserving a certain
pinching condition in dimension four was proved by Hamilton [H
5,\S 4]; his argument works in our case too (and the estimates are
much easier to verify). The point is that we can change our
$\de$-neck metric near the middle of the neck by a conformal
factor $e^{-f},$ where $f=f(z)$ is positive on the part of the
neck we want to remove, and zero on the part we want to preserve,
and $z$ is the coordinate along $\mathbb{I}$ in our
parametrization $\mathbb{S}^2\times\mathbb{I}$ of the neck. Then,
in the region near the middle of the neck, where $f$ is small, the
dominating terms in the formulas for the change of curvature are
just positive constant multiples of $f'',$ so the pinching
improves, and  all the curvatures become positive on the set where
$f>\de'.$
\par Now we can continue our solution until it becomes singular
for the next time. Note that after the surgery the manifold may
become disconnected; in this case, each component should be dealt
with separately. Furthermore, let us agree to declare extinct
every component which is $\e$-close to a metric quotient of the
round sphere; that allows to exclude such components from the list
of canonical neighborhoods. Now since every surgery reduces the
volume by at least $h^3,$ the sequence of surgery times is
discrete, and, taking for granted the a priori assumptions, we can
continue our solution indefinitely, not ruling out the possibility
that it may become extinct at some finite time.
\par {\bf 4.5} In order to justify the canonical neighborhood
assumption in the next section, we need to check several
assertions.
\par {\bf Lemma.} {\it For any $A<\infty ,0< \theta<1, $ one
can find $\bar{\de}=\bar{\de}(A,\theta)$ with the following
property. Suppose we have a solution to the Ricci flow with
$\de$-cutoff, satisfying the a priori assumptions on $[0,T],$ with
$\de<\bar{\de}.$ Suppose we have a surgery at time $T_0\in (0,T),$
let $p$ correspond to the center of the standard cap, and let
$T_1=\mathrm{min}(T,T_0+\theta h^2).$ Then either
\par (a) The solution is defined on $P(p,T_0,Ah,T_1-T_0),$ and is,
after scaling with factor $h^{-2}$ and shifting time $T_0$ to
zero, $A^{-1}$-close to the corresponding subset on the standard
solution from section 2, or
\par (b) The assertion (a) holds with $T_1$ replaced by some time
$t^+\in [T_0,T_1),$ where $t^+$ is a surgery time; moreover, for
each point in $B(p,T_0,Ah),$ the solution is defined for $t\in
[T_0,t^+)$ and is not defined past $t^+.$}
\par {\it Proof.} Let $Q$ be the maximum of the scalar curvature
on the standard solution in the time interval $[0,\theta],$ let
$\triangle t=N^{-1}(T_1-T_0)<\e\eta^{-1}Q^{-1}h^2,$ and let
$t_k=T_0+k\triangle t, k=0,...,N.$
\par Assume first that for each
point in $B(p,T_0,A_0h),$ where $A_0=\e(\de ')^{-1},$ the solution
is defined on $[t_0,t_1].$ Then by (1.3) and the choice of
$\triangle t$ we have a uniform curvature bound on this set for
$h^{-2}$-scaled metric. Therefore we can define  $A_1,$ depending
only on $A_0$ and tending to infinity with $A_0,$ such that the
solution in $P(p,T_0,A_1h,t_1-t_0)$ is, after scaling and time
shifting, $A_1^{-1}$-close to the corresponding subset in the
standard solution. In particular, the scalar curvature on this
subset does not exceed $2Qh^{-2}.$ Now if for each point in
$B(p,T_0,A_1h)$ the solution is defined on $[t_1,t_2],$ then we
can repeat the procedure, defining $A_2$ etc. Continuing this way,
we eventually define $A_N,$ and it would remain to choose $\de$ so
small, and correspondingly $A_0$ so large, that $A_N>A.$
\par Now assume that for some $k, 0\le k<N,$ and for some $x\in
B(p,T_0,A_kh)$ the solution is defined on $[t_0,t_k]$ but not on
$[t_k,t_{k+1}].$ Then we can find a surgery time
$t^+\in[t_k,t_{k+1}],$ such that the solution on $B(p,T_0,A_kh)$
is defined on $[t_0,t^+),$ but for some points of this ball it is
not defined past $t^+.$ Clearly, the $A_{k+1}^{-1}$-closeness
assertion holds on $P(p,T_0,A_{k+1}h,t^+-T_0).$ On the other hand,
the solution on $B(p,T_0,A_kh)$ is at least $\e$-close to the
standard one for all $t\in [t_k,t^+),$ hence no point of this set
can be the center of a $\de$-neck neighborhood at time $t^+.$
However, the surgery is always done along the middle two-sphere of
such a neck. It follows that for each point of $B(p,T_0,A_kh)$ the
solution terminates at $t^+.$
\par {\bf 4.6 Corollary.} {\it For any $l<\infty$ one can find
$A=A(l)<\infty$ and $\theta=\theta(l), 0<\theta<1,$ with the
following property. Suppose we are in the situation of the lemma
above, with $\de<\bar{\de}(A,\theta).$ Consider smooth curves
$\gamma$ in the set $B(p,T_0,Ah),$ parametrized by $t\in
[T_0,T_{\gamma}],$ such that $\gamma(T_0)\in B(p,T_0,Ah/2)$ and
either $T_{\gamma}=T_1<T$, or $T_{\gamma}<T_1$ and
$\gamma(T_{\gamma})\in\partial B(p,T_0,Ah).$ Then
$\int_{T_0}^{T_{\gamma}}{(R(\gamma(t),t)+|\dot{\gamma}(t)|^2)dt}>l$.}
\par {\it Proof.} Indeed, if $T_{\gamma}=T_1,$ then on the
standard solution we would have
$\int_{T_0}^{T_{\gamma}}R(\gamma(t),t)dt
\ge\co\int_0^{\theta}(1-t)^{-1}dt=-\co\cdot
(\mathrm{log}(1-\theta))^{-1},$ so by choosing $\theta$
sufficiently close to one we can handle this case. Then we can
choose $A$
 so large that on the standard solution
 $\mathrm{dist}_t(p,\partial B(p,0,A))\ge 3A/4$ for each $t\in
 [0,\theta].$ Now if $\gamma(T_{\gamma})\in\partial B(p,T_0,Ah)$
 then $\int_{T_0}^{T_{\gamma}}|\dot{\gamma}(t)|^2dt\ge A^2/100,$
 so by taking $A$ large enough, we can handle this case as well.
 \par {\bf 4.7 Corollary.} {\it For any $Q<\infty$ there exists
 $\theta=\theta(Q), 0<\theta<1$ with the following property.
 Suppose we are in the situation of the lemma above, with
 $\de<\bar{\de}(A,\theta), A>\e^{-1}.$   Suppose that for some point $x\in
 B(p,T_0,Ah)$ the solution is defined at $x$ (at least) on $[T_0,T_x], T_x\le T,$
 and satisfies $Q^{-1}R(x,t)\le R(x,T_x)\le Q(T_x-T_0)^{-1}$ for
 all $t\in [T_0,T_x].$
 Then $T_x\le T_0 + \theta h^2.$}
 \par {\it Proof.} Indeed, if $\ \ \ T_x>T_0+\theta h^2,$ then by lemma
 $\ \ \ R(x,T_0+\theta h^2)\ge \co\cdot (1-\theta)^{-1}h^{-2},$ whence
 $R(x,T_x)\ge \co\cdot Q^{-1}(1-\theta)^{-1}h^{-2},$ and $T_x-T_0\le
 \co\cdot Q^2(1-\theta)h^2<\theta h^2$ if $\theta$ is close enough to
 one.
 \section{Justification of the a priori assumption}
 \par {\bf 5.1} Let us call a riemannian manifold $(M,g_{ij})$
 normalized if $M$ is a closed oriented 3-manifold, the sectional
 curvatures of $g_{ij}$ do not exceed one in absolute value, and
 the volume of every metric ball of radius one is at least half
 the volume of the euclidean unit ball. For smooth Ricci flow with
 normalized initial data we have, by [H 4, 4.1], at any time $t>0$ the pinching
 estimate \begin{equation}
 Rm\ge-\phi(R(t+1))R, \end{equation} where $\phi$ is a decreasing function, which
 behaves at infinity like $\frac{1}{\mathrm{log}}.$
 As explained in 4.4, this pinching estimate can be preserved for
 Ricci flow with $\de$-cutoff. Justification of the canonical
 neighborhood assumption requires additional arguments.
 In fact, we are able to construct solutions satisfying this
 assumption only allowing $r$ and $\de$ be
 functions of time rather than constants; clearly, the arguments
 of the previous section are valid in this case, if we assume that
 $r(t),\ \de(t)$ are non-increasing, and bounded
 away from zero on every finite time interval.
 \par {\bf Proposition.} {\it There exist decreasing sequences $0<r_j<\e^2, \ka_j>0,
 0<\bar{\de}_j<\e^2, j=1,2,...,$ such that for any
 normalized initial data and any function $\de(t),$ satisfying
 $0<\de(t)<\bar{\de}_j$ for $t\in [2^{j-1}\e,2^j\e],$ the Ricci flow
 with $\de(t)$-cutoff is defined for $t\in[0,+\infty]$ and satisfies
 the $\ka_j$-noncollapsing assumption and the
 canonical neighborhood assumption with parameter $r_j$
 on the time interval $[2^{j-1}\e,2^j\e].$( Recall that we have excluded from
 the list of canonical neighborhoods the closed manifolds, $\e$-close to metric
 quotients of the round sphere. Complete extinction of the
 solution
 in finite time is not ruled out.)}
 \par The proof of the proposition is by induction: having constructed our sequences
 for $1\le j\le i,$ we make one more step, defining
 $r_{i+1}, \ka_{i+1}, \bar{\de}_{i+1},$ and redefining $\bar{\de}_i=\bar{\de}_{i+1};$
 each step is  analogous to the proof of
 Theorem I.12.1.
 \par First we need to check a $\kappa$-noncollapsing
 condition.
 \par {\bf 5.2 Lemma.} {\it Suppose we have constructed the sequences,
 satisfying the proposition for $1\le j\le i.$ Then there exists $\ka>0,$
 such that for any $r, 0<r<\e^2,$ one can find
 $\bar{\de}=\bar{\de}(r)>0,$ which may also depend on the already constructed
 sequences, with the following property. Suppose we have a
 solution to the Ricci flow with $\de(t)$-cutoff on a time interval
 $[0,T],$ with normalized initial data, satisfying the proposition
 on $[0,2^i\e],$ and the canonical neighborhood assumption with
 parameter $r$ on $[2^i\e,T],$ where $2^i\e\le T\le 2^{i+1}\e,\ \
  0<\de(t)<\bar{\de}\ \  for \ \  t\in [2^{i-1}\e,T] .$ Then it
  is
 $\kappa$-noncollapsed on all scales less than $\e.$}
 \par {\it Proof.} Consider a neighborhood
 $P(x_0,t_0,r_0,-r_0^2), 2^i\e<t_0\le T, 0<r_0<\e, $ where the solution is defined and
 satisfies $|Rm|\le r_0^{-2}.$ We may assume $r_0\ge r,$
 since otherwise the lower bound for the volume of the ball
 $B(x_0,t_0,r_0) $ follows from the canonical neighborhood
 assumption. If the solution was smooth
 everywhere, we could estimate from below the volume of the ball
 $B(x_0,t_0,r_0)$ using the argument from [I.7.3]: define
 $\tau(t)=t_0-t$ and consider the reduced volume function using
 the $\mathcal{L}$-exponential map from $x_0;$ take a point
 $(x,\e)$ where the reduced distance $l$ attains its minimum for $\tau=t_0-\e,$
 $\ \ l(x,\tau)\le 3/2;$ use it to obtain an upper bound for the
 reduced distance to the points of $B(x,0,1),$ thus getting a
 lower bound for the reduced volume at $\tau=t_0,$ and apply the
 monotonicity formula. Now if the solution undergoes surgeries,
 then we still can measure the $\mathcal{L}$-length, but only for
 admissible curves, which stay in the region,  unaffected by
 surgery. An inspection of the constructions in [I,\S 7] shows
 that the argument would go through if we knew that every barely
 admissible curve, that is a curve on the boundary of the set of
 admissible curves, has reduced length at least $3/2+\kappa '$ for
 some fixed $\kappa '>0.$ Unfortunately, at the moment I don't see
  how to ensure that without imposing new restrictions on
   $\de(t)$ for all $t\in [0,T],$ so we need some additional arguments.
 \par Recall that for a curve $\gamma,$ parametrized by $t,$ with
 $\gamma(t_0)=x_0,$ we have $\mathcal{L}(\gamma ,\tau)=
 \int_{t_0-\tau}^{t_0}{\sqrt{t_0-t}(R(\gamma(t),t)+|\dot{\gamma}(t)|^2)dt}.$
 We can also define $\mathcal{L}_+(\gamma,\tau)$ by replacing in the previous formula $R$
 with $R_+=\mathrm{max}(R,0).$ Then $\mathcal{L}_+\le
 \mathcal{L}+4T\sqrt{T}$ because $R\ge -6$ by the maximum
 principle and normalization. Now suppose we could show that every
 barely admissible curve with endpoints $(x_0,t_0)$ and $(x,t),$
 where $t\in [2^{i-1}\e,T),$ has $\mathcal{L}_+>
 2\e^{-2}T\sqrt{T};$ then we could argue that either there exists a
 point $(x,t), t\in [2^{i-1}\e,2^i\e],$ such that $R(x,t)\le
 r_i^{-2}$ and $\mathcal{L}_+\le \e^{-2}T\sqrt{T},$ in which case
 we can take this point in place of $(x,\e)$ in the argument of
 the previous paragraph, and obtain (using Claim 1 in 4.2) an estimate for $\ka$ in terms
 of $r_i,\ka_i,T,$ or for any $\gamma,$ defined on
 $[2^{i-1}\e,t_0], \gamma(t_0)=x_0,$ we have $\mathcal{L}_+\ge
\mi(\e^{-2}T\sqrt{T}, \frac{2}{3}(2^{i-1}\e)^{\frac{3}{2}}
r_i^{-2})>\e^{-2}T\sqrt{T},$ which is in contradiction with the
assumed bound for barely admissible curves and the bound $\mi\
l(x,t_0-2^{i-1}\e)\le 3/2,$ valid in the smooth case.
 Thus,  to conclude the proof it is
 sufficient to check the following assertion.
 \par {\bf 5.3 Lemma.} {\it For any $\mathcal{L}<\infty$ one can
 find $\bar{\de}=\bar{\de}(\mathcal{L},r_0)>0$ with the
 following property. Suppose that in the situation of the previous
 lemma we have a curve $\gamma,$
 parametrized by $t\in[T_0,t_0], 2^{i-1}\e\le T_0<t_0,$ such that $\gamma(t_0)=x_0,$
 $T_0$ is a surgery time, and $\gamma(T_0)\in B(p,T_0,\e^{-1}h),$
 where $p$ corresponds to the center of the cap, and $h$ is the
 radius of the $\de$-neck. Then we have an estimate
 $\int_{T_0}^{t_0}{\sqrt{t_0-t}(R_+(\gamma(t),t)+|\dot{\gamma}(t)|^2)dt}\ge\mathcal{L}.$}
 \par {\it Proof.} It is clear that if we take $\triangle t=\e
 r_0^4\mathcal{L}^{-2},$ then either $\gamma$ satisfies our
 estimate, or $\gamma$ stays in $P(x_0,t_0,r_0,-\triangle t)$ for
 $t\in [t_0-\triangle t,t_0].$ In the latter case  our estimate
 follows from Corollary 4.6, for $l=\mathcal{L}(\triangle
 t)^{-\frac{1}{2}},$ since clearly $T_{\gamma}<t_0-\triangle t$
 when $\de$ is small enough.
 \par {\bf 5.4} {\it Proof of proposition.} Assume the contrary,
 and let the sequences $r^{\al}, \bar{\de}^{\al\be}$
  be such that $r^{\al}\to 0$ as
 $\al\to\infty,$ $\ \ \bar{\de}^{\al\be}\to 0$ as $\be\to\infty$ with
 fixed $\al,$ and let $(M^{\al\be},g_{ij}^{\al\be})$ be normalized
 initial data for solutions to the Ricci flow with $\de(t)$-cutoff,
 $\de(t)<\bar{\de}^{\al\be}$
 on $[2^{i-1}\e,2^{i+1}\e],$ which satisfy
 the statement on $[0,2^i\e],$ but
 violate the canonical neighborhood assumption with parameter $r^{\al}$
 on $[2^i\e,2^{i+1}\e].$ Slightly abusing notation,
we'll drop the indices $\al ,\be$ when we consider an individual
solution.
 \par  Let
 $\bar{t}$ be the first time when the assumption is violated
 at some point $\bar{x};$ clearly
 such time exists, because it is an open condition. Then by lemma 5.2 we have
 uniform $\ka$-noncollapsing on $[0,\bar{t}].$ Claims
 1,2 in 4.2
 are also valid on $[0,\bar{t}];$ moreover, since $h<<r,$ it follows from Claim 1 that
 the solution is defined on the whole parabolic neighborhood indicated there in case
 $R(x_0,t_0)\le r^{-2}.$
 \par Scale
 our solution about $(\bar{x},\bar{t})$ with factor $R(\bar{x},\bar{t})\ge r^{-2}$ and take a
 limit for subsequences of $\al,\be\to\infty.$  At time $\bar{t},$ which we'll shift
 to zero in the limit, the curvature bounds
 at finite distances from $\bar{x}$ for the scaled metric are ensured
 by Claim 2 in 4.2. Thus, we get  a smooth complete limit of
 nonnegative sectional curvature, at time zero. Moreover, the
 curvature of the limit is uniformly bounded, since otherwise it
 would contain $\e$-necks of arbitrarily small radius.
 \par Let $Q_0$
 denote the curvature bound. Then, if there was no surgery, we
 could, using Claim 1 in 4.2, take a limit on the time interval
 $[-\e\eta^{-1}Q_0^{-1},0].$ To prevent this, there must exist
 surgery times $T_0\in
 [\bar{t}-\e\eta^{-1}Q_0^{-1}R^{-1}(\bar{x},\bar{t}),\bar{t}]$ and
 points $x$ with
 $\mathrm{dist}^2_{T_0}(x,\bar{x})R^{-1}(\bar{x},\bar{t}) $
 uniformly bounded as $\al,\be\to\infty ,$ such that the solution at $x$ is defined on
 $[T_0,\bar{t}],$ but not before $T_0.$ Using Claim 2 from 4.2 at time $T_0,$
 we see that $R(\bar{x},\bar{t})h^2(T_0)$ must be bounded away from zero.
 Therefore, in this case we can apply
 Corollary 4.7, Lemma 4.5 and Claim 5 in section 2 to show that the point
 $(\bar{x},\bar{t})$ in
 fact has a canonical neighborhood, contradicting its choice.
 (It is not excluded that the strong $\e$-neck neighborhood extends to times before
 $T_0,$ where it is a part of the strong $\de$-neck that existed
 before surgery.)
 \par Thus we have a limit on a certain time interval. Let $Q_1$
 be the curvature bound   for   this   limit.   Then we   either   can
 construct a limit on the time interval
 $[-\e\eta^{-1}(Q_0^{-1}+Q_1^{-1}),0],$ or there
 is a surgery, and we get a contradiction as before. We can
 continue this procedure indefinitely, and the final part of the
 proof of Theorem I.12.1 shows that the  bounds $Q_k$ can not go to
 infinity while the limit is defined on a bounded time interval. Thus we get a limit on
 $(-\infty,0],$ which is $\ka$-noncollapsed by Lemma 5.2, and
 this means that $(\bar{x},\bar{t})$ has a canonical neighborhood
 by the results of section 1 - a contradiction.
 \section{Long time behavior I}
 \par {\bf 6.1} Let us summarize what we have achieved so far.
 We have shown the existence of decreasing (piecewise constant) positive functions
 $r(t)$ and $\bar{\de}(t)$ (which we may assume converging
 to zero at infinity), such that if $ (M,g_{ij})$ is a
 normalized manifold, and $0<\de(t)<\bar{\de}(t),$ then there exists
 a solution to the Ricci flow with $\de(t)$-cutoff on the time interval $[0,+\infty],$
 starting from
 $(M,g_{ij})$ and satisfying on each subinterval $[0,t]$ the
 canonical neighborhood assumption with parameter $r(t),$ as well
 as the pinching estimate (5.1).
 \par In particular, if the initial
 data has positive scalar curvature, say $R\ge a>0,$ then the
 solution becomes extinct in time at most $\frac{3}{2a},$ and it
 follows that $M$ in this case is diffeomorphic to a connected sum
 of several copies of $\mathbb{S}^2\times\mathbb{S}^1$ and metric
 quotients of round $\mathbb{S}^3.$ ( The topological description of 3-manifolds with
 positive scalar curvature modulo quotients of homotopy spheres
 was obtained by Schoen-Yau and Gromov-Lawson more than 20 years ago, see [G-L] for
 instance; in particular, it is well known and
 easy to check that every  manifold that can be decomposed in a connected sum
 above admits a metric of positive scalar
 curvature.) Moreover, if the scalar
 curvature is only nonnegative, then by the strong maximum
 principle it instantly becomes positive unless the metric is
 (Ricci-)flat; thus in this case, we need to add to our list the
 flat manifolds.
 \par However, if the scalar curvature is negative
 somewhere, then we need to work more in order to understand the
 long tome behavior of the solution. To achieve this we need first
 to prove versions of Theorems I.12.2 and I.12.3 for solutions
 with cutoff.
 \par {\bf 6.2 Correction to Theorem I.12.2.} Unfortunately, the
 statement of Theorem I.12.2 was incorrect. The assertion I had in
 mind is as follows:
 \par {\it Given a function $\phi$ as above, for any $A<\infty$
 there exist $K=K(A)<\infty$ and $\rho=\rho(A)>0$ with the
 following property. Suppose in dimension three we have a solution
 to the Ricci flow with $\phi$-almost nonnegative curvature, which
 satisfies the assumptions of theorem 8.2 for some $x_0,r_0$ with
 $\phi(r_0^{-2})<\rho.$
 Then $R(x,r_0^2)\le Kr_0^{-2}$ whenever
 $\mathrm{dist}_{r_0^2}(x,x_0)<Ar_0.$}
 \par It is this assertion that was used in the proof of Theorem I.12.3
 and Corollary I.12.4.
 \par {\bf 6.3 Proposition.} {\it For any $A<\infty$ one can find
 $\ka=\ka(A)>0, K_1=K_1(A)<\infty ,K_2=K_2(A)<\infty ,
 \bar{r}=\bar{r}(A)>0,$ such that for any $t_0<\infty$ there exists
 $\bar{\de}=\bar{\de}_A(t_0)>0,$ decreasing in $t_0, $ with the following property. Suppose
 we have a solution to the Ricci flow with $\de(t)$-cutoff on time interval
 $[0,T],\ \
 \de(t)<\bar{\de}(t)$ on $[0,T],\  \de(t)<\bar{\de}$ on $[t_0/2,t_0] , $ with
 normalized initial data; assume that the solution is defined
 in the whole parabolic neighborhood $P(x_0,t_0,r_0,-r_0^2), \ 2r_0^2<t_0,$ and
 satisfies $|Rm|\le r_0^{-2}$ there, and that the volume of the
 ball $B(x_0,t_0,r_0)$ is at least $A^{-1}r_0^3.$ Then
 \par (a) The solution is $\ka$-noncollapsed on the scales less
 than $r_0$ in the ball $B(x_0,t_0,Ar_0).$
 \par (b) Every point $x\in
 B(x_0,t_0,Ar_0)$ with $R(x,t_0)\ge K_1 r_0^{-2}$ has a canonical
 neighborhood as in 4.1.
 \par (c) If $r_0\le \bar{r}\sqrt{t_0}$ then $R\le K_2r_0^{-2}$ in
 $B(x_0,t_0,Ar_0).$}
\par {\it Proof.} (a) This is an analog of Theorem I.8.2.
Clearly we have $\ka$-noncollapsing on the scales less than
$r(t_0),$ so we may assume $r(t_0)\le r_0\le \sqrt{t_0/2}$ , and
study the scales $\rho, r(t_0)\le \rho\le r_0.$ In particular, for
fixed $t_0$ we are interested in the scales, uniformly equivalent
to one. \par  So assume that $x\in B(x_0,t_0,Ar_0)$ and the
solution is defined in the whole $P(x,t_0,\rho,-\rho^2)$ and
satisfies $|Rm|\le \rho^{-2}$ there. An inspection of the proof of
I.8.2 shows that in order to make the argument work it suffices to
check that for any barely admissible curve $\gamma$, parametrized
by $t\in [t_{\gamma},t_0], t_0-r_0^2\le t_{\gamma}\le t_0,$ such
that $\gamma(t_0)=x,$ we have an estimate \begin{equation}
2\sqrt{t_0-t_{\gamma}}
\int_{t_{\gamma}}^{t_0}{\sqrt{t_0-t}(R(\gamma(t),t)+|\dot{\gamma}(t)|^2)dt}\ge
C(A)r_0^2 \end{equation} for a certain function $C(A)$ that can be
made explicit. Now we would like to conclude the proof by using
Lemma 5.3. However, unlike the situation in Lemma 5.2, here Lemma
5.3 provides the estimate we need only if $t_0-t_{\gamma}$ is
bounded away from zero, and otherwise we only get an estimate $
\rho^2$ in place of $C(A)r_0^2.$ Therefore we have to return to
the proof of I.8.2.
\par Recall that in that proof we scaled the solution to make
$r_0=1$ and worked on the time interval $[1/2,1].$ The maximum
principle for the evolution equation of the scalar curvature
implies that on this time interval we have $R\ge -3.$ We
considered a function of the form
$h(y,t)=\phi(\hat{d}(y,t))\hat{L}(y,\tau),$ where $\phi$ is a
certain cutoff function, $\tau=1-t,
\hat{d}(y,t)=\mathrm{dist}_t(x_0,y)-A(2t-1),
\hat{L}(y,\tau)=\bar{L}(y,\tau)+7,$ and $\bar{L}$ was defined in
[I,(7.15)]. Now we redefine $\hat{L},$ taking
$\hat{L}(y,\tau)=\bar{L}(y,\tau)+2\sqrt{\tau}.$  Clearly,
$\hat{L}>0$ because $R\ge -3$ and $2\sqrt{\tau}>4\tau^2$ for
$0<\tau\le 1/2.$ Then the computations and estimates of I.8.2
yield $$ \Box h\ge -C(A)h-(6+\frac{1}{\sqrt{\tau}})\phi$$ Now
denoting by $h_0(\tau)$ the minimum of $h(y,1-t),$ we can estimate
\begin{equation}
 \frac{d}{d\tau}(\mathrm{log}(\frac{h_0(\tau)}{\sqrt{\tau}}))\le
C(A)+\frac{6\sqrt{\tau}+1}{2\tau-4\tau^2\sqrt{\tau}}-\frac{1}{2\tau}\le
C(A)+ \frac{50}{\sqrt{\tau}}, \end{equation} whence
\begin{equation} h_0(\tau)\le \sqrt{\tau}\ \mathrm{exp}(C(A)\tau
+100\sqrt{\tau}),\end{equation}  because the left hand side of
(6.2)  tends to zero as $\tau\to 0+.$
\par Now we can return to our proof, replace the right hand
side of (6.1) by the right hand side of (6.3) times $r_0^2,$ with
$\tau=r_0^{-2}(t_0-t_{\gamma}),$ and apply Lemma 5.3.
 \par (b) Assume the contrary, take a sequence
 $K_1^{\al}\to\infty$ and consider the solutions violating the
 statement. Clearly, $K_1^{\al} (r_0^{\al})^{-2}<(r(t_0^{\al}))^{-2},$ whence
 $t_0^{\al}\to\infty;$
\par When $K_1$ is
large enough, we can,
  arguing as in the proof
 of Claim 1 in [I.10.1], find a point $(\bar{x},\bar{t}),
 x\in B(x_0,\bar{t},2Ar_0), \bar{t}\in [t_0-r_0^2/2,t_0],$ such
 that $\bar{Q}=R(\bar{x},\bar{t})>K_1r_0^{-2},\ \  (\bar{x},\bar{t})$ does
 not satisfy the canonical neighborhood assumption, but each point
 $(x,t)\in \bar{P}$ with $R(x,t)\ge 4\bar{Q}$ does, where $\bar{P}$
 is the set of all $(x,t)$ satisfying $\bar{t}-\frac{1}{4}K_1\bar{Q}^{-1}\le t\le
 \bar{t},\ \
 \mathrm{dist}_{t}(x_0,x)\le\mathrm{dist}_{\bar{t}}(x_0,\bar{x})+K_1^{\frac{1}{2}}
 \bar{Q}^{-\frac{1}{2}}.$ (Note that $\bar{P}$ is not a parabolic
 neighborhood.) Clearly we can use (a) with slightly different
 parameters to ensure $\ka$-noncollapsing in $\bar{P}.$
 \par Now we apply the argument from 5.4. First, by Claim 2 in 4.2,
 for any $\bar{A}<\infty$ we have an estimate $R\le
 Q(\bar{A})\bar{Q}$ in
 $B(\bar{x},\bar{t},\bar{A}\bar{Q}^{-\frac{1}{2}})$ when $K_1$ is large enough; therefore we
 can take a limit as $\al\to\infty$ of scalings with factor $\bar{Q}$
 about $(\bar{x},\bar{t}),$ shifting the time $\bar{t}$ to zero;
 the limit at time zero would be a smooth complete nonnegatively curved
 manifold. Next we observe that this limit has curvature uniformly
 bounded, say, by $Q_0,$ and therefore, for each fixed $\bar{A}$
 and for sufficiently large $K_1,$ the parabolic neighborhood
 $P(\bar{x},\bar{t},\bar{A}\bar{Q}^{-\frac{1}{2}},-\e\eta^{-1}Q_0^{-1}\bar{Q}^{-1})$
 is contained in $\bar{P}.$  (Here we use the estimate of distance change,
 given by Lemma I.8.3(a).) Thus we can take a limit on the
 interval $[-\e\eta^{-1}Q_0^{-1},0]. $ (The possibility of
 surgeries is ruled out as in 5.4) Then we repeat the procedure
 indefinitely, getting an ancient $\ka$-solution in the limit,
 which means a contradiction.
 \par (c) If $x\in B(x_0,t_0,Ar_0)$ has very large curvature, then
  on the shortest geodesic $\gamma$ at time
 $t_0,$ that connects $x_0$ and $x,$ we can find a point $y,$ such
 that $R(y,t_0)=K_1(A)r_0^{-2}$ and the curvature is larger at all
 points of the segment of $\gamma$ between $x$ and $y.$ Then our
 statement follows from Claim 2 in 4.2, applied to this segment.
 \par From now on we redefine the function $\bar{\de}(t)$ to be
 $\mathrm{min}(\bar{\de}(t),\bar{\de}_{2t}(2t)),$   so  that the
 proposition above always holds for $A=t_0.$
 \par {\bf 6.4 Proposition.} {\it There exist $\tau>0, \bar{r}>0,
 K<\infty$ with the following property. Suppose we have a solution
 to the Ricci flow with $\de(t)$-cutoff on the time interval $[0,t_0],$
  with normalized initial data. Let $r_0,t_0$ satisfy
  $2C_1h\le r_0\le \bar{r}\sqrt{t_0},$ where $h$ is the maximal cutoff radius
  for surgeries in $[t_0/2,t_0],$ and assume that
   the ball $B(x_0,t_0,r_0)$ has sectional curvatures at least
   $-r_0^{-2}$ at each point, and the volume of any subball
   $B(x,t_0,r)\subset B(x_0,t_0,r_0)$ with any radius $r>0$ is at
   least $(1-\e)$ times the volume of the euclidean ball of the
   same radius. Then the solution is defined in
   $P(x_0,t_0,r_0/4,-\tau r_0^2)$ and satisfies $R< Kr_0^{-2}$ there.}
   \par {\it Proof.} Let us first consider the case $r_0\le
   r(t_0).$ Then clearly $R(x_0,t_0)\le C_1^2r_0^{-2},$ since an
   $\e$-neck of radius $r$ can not contain an almost euclidean
   ball of radius $\ge r.$ Thus we can take $K=2C_1^2,
   \tau=\e\eta^{-1}C_1^{-2}$ in this case, and since $r_0\ge 2C_1h,$
   the surgeries do not interfere in
   $P(x_0,t_0,r_0/4,-\tau r_0^2).$
   \par In order to handle the other case
   $r(t_0)<r_0\le\bar{r}\sqrt{t_0}$ we need a couple of lemmas.
   \par {\bf 6.5 Lemma.} {\it There exist $\tau_0>0$ and $K_0<\infty,$
   such that if we have a smooth solution to the Ricci flow in
   $P(x_0,0,1,-\tau),\tau\le\tau_0,$ having sectional curvatures
   at least $-1$, and the volume of the ball $B(x_0,0,1)$ is at
   least $(1-\e)$ times the volume of the euclidean unit ball,
   then
   \par (a) $R\le K_0\tau^{-1}$ in $P(x_0,0,1/4,-\tau/2),$ and
   \par (b) the ball $B(x_0,1/4,-\tau)$ has volume at least
   $\frac{1}{10}$ times the volume of the euclidean ball of the
   same radius.}
   \par The proof can be extracted from the proof of Lemma I.11.6.
   \par {\bf 6.6 Lemma.} {\it For any $w>0$ there exists
   $\theta_0=\theta_0(w)>0,$ such that if $B(x,1)$ is a metric
   ball of volume at least $w,$ compactly contained in a manifold
   without boundary with sectional curvatures at least $-1,$ then
   there exists a ball $B(y,\theta_0)\subset B(x,1),$ such that
   every subball $B(z,r)\subset B(y,\theta_0)$ of any radius $r$
   has volume at least $(1-\e)$ times the volume of the euclidean
   ball of the same radius.}
   \par This is an elementary fact from the theory of Aleksandrov
   spaces.
   \par {\bf 6.7} Now we continue the proof of the proposition. We claim
   that one can take $\tau=\mathrm{min}(\tau_0/2,
   \e\eta^{-1}C_1^{-2}), K=\mathrm{max}(2K_0\tau^{-1}, 2C_1^2).$
   Indeed, assume the contrary, and take a sequence of
   $\bar{r}^{\al}\to 0$ and solutions, violating our assertion for
   the chosen $\tau,K.$ Let $t_0^{\al}$ be the first time when it
   is violated, and let $B(x_0^{\al},t_0^{\al},r_0^{\al})$ be the
    counterexample with the smallest radius. Clearly
    $r_0^{\al}>r(t_0^{\al})$ and $(r_0^{\al})^2 (t_0^{\al})^{-1}\to
    0$ as $\al\to\infty.$
    \par Consider any ball $B(x_1,t_0,r)\subset B(x_0,t_0,r_0), r<r_0.$
    Clearly we can apply our proposition to this ball and get the
    solution in $P(x_1,t_0,r/4,-\tau r^2)$ with the curvature
    bound $R<Kr^{-2}.$ Now if $r_0^2 t_0^{-1}$ is small enough,
    then we can apply proposition 6.3(c) to get an estimate $R(x,t)\le
    K'(A)r^{-2}$ for $(x,t)$ satisfying $t\in [t_0-\tau
    r^2/2,t_0], \mathrm{dist}_t(x,x_1)<Ar,$ for some function
    $K'(A)$ that can be made explicit. Let us choose
    $A=100r_0r^{-1};$ then we get the solution with a curvature
    estimate in $P(x_0,t_0,r_0,-\triangle t),$ where $\triangle
    t=K'(A)^{-1}r^2.$ Now the pinching estimate implies $Rm\ge
    -r_0^{-2}$ on this set, if $r_0^2t_0^{-1}$ is small enough
    while $rr_0^{-1}$ is bounded away from zero. Thus we can use
    lemma 6.5(b) to estimate the volume of the ball
    $B(x_0,t_0-\triangle t,r_0/4)$ by at least $\frac{1}{10}$ of
    the volume of
    the  euclidean ball of the same radius, and then by lemma 6.6
    we can find a subball $B(x_2,t_0-\triangle t,
    \theta_0(\frac{1}{10})r_0/4),$ satisfying the assumptions of
    our proposition. Therefore, if we put
    $r=\theta_0(\frac{1}{10})r_0/4,$ then we can repeat our
    procedure as many times as we like, until we reach the time
    $t_0-\tau_0 r_0^2$, when the lemma 6.5(b) stops working. But
    once we reach this time, we can apply lemma 6.5(a) and get the
    required curvature estimate, which is a contradiction.
    \par {\bf 6.8 Corollary.} {\it For any $w>0$ one can find
    $\tau=\tau(w)>0, K=K(w)<\infty, \bar{r}=\bar{r}(w)>0, \theta=\theta(w)>0$ with the
    following property. Suppose we have a solution to the Ricci
    flow with $\de(t)$-cutoff on the time interval $[0,t_0],$ with
    normalized initial data. Let $t_0,r_0$ satisfy
    $\theta^{-1}(w)h\le r_0\le \bar{r}\sqrt{t_0},$ and assume that
    the ball $B(x_0,t_0,r_0)$ has sectional curvatures at least
    $-r_0^2$ at each point, and volume at least $wr_0^3.$ Then
    the solution is defined in $P(x_0,t_0,r_0/4,-\tau r_0^2)$ and
    satisfies $R<Kr_0^{-2}$ there.}
    \par Indeed, we can apply proposition 6.4 to a smaller ball,
    provided by lemma 6.6, and then use proposition 6.3(c).
    \section{Long time behavior II}
    \par In this section we adapt the arguments of Hamilton [H 4]
    to a more general setting. Hamilton considered smooth Ricci flow
    with bounded normalized curvature; we drop both these assumptions.
    In the end of [I,13.2] I claimed that the volumes of the maximal horns
    can be effectively bounded below, which would imply that the solution must be
    smooth from some time on; however, the argument I had in mind
    seems to be faulty. On the other hand, as we'll see below, the presence of
     surgeries does not lead to any substantial problems.
     \par   From now on we assume that our
    initial manifold does not admit a metric with nonnegative
    scalar curvature, and that once we get a component with
    nonnegative scalar curvature, it is immediately removed.
    \par {\bf 7.1} (cf. [H 4,\S 2,7]) Recall that for a solution to the smooth Ricci flow the scalar
    curvature satisfies the evolution equation \begin{equation}
    \frac{d}{dt}R=\triangle R+2|Ric|^2= \triangle R
    +2|Ric^{\circ}|^2+\frac{2}{3}R^2, \end{equation}
    where $Ric^{\circ}$ is the trace-free part of $Ric.$ Then
    $R_{\mi}(t)$ satisfies
    $\frac{d}{dt}R_{\mi}\ge \frac{2}{3}R_{\mi}^2,$
    whence \begin{equation} R_{\mi}(t)\ge
    -\frac{3}{2}\ \ \frac{1}{t+1/4} \end{equation} for a solution with
    normalized initial data. The evolution equation for the volume
    is $\frac{d}{dt}V=-\int RdV,$
 in particular \begin{equation} \frac{d}{dt}V\le -R_{\mi}V,
 \end{equation} whence by (7.2) the function
 $V(t)(t+1/4)^{-\frac{3}{2}}$ is non-increasing in $t.$ Let
 $\bar{V}$ denote its limit as $t\to\infty.$
 \par Now the scale invariant quantity
 $\hat{R}=R_{\mi}V^{\frac{2}{3}}$ satisfies \begin{equation}
 \frac{d}{dt}\hat{R}(t)\ge \frac{2}{3}\ \ \hat{R}V^{-1}\int
 (R_{\mi}-R)dV , \end{equation} which is nonnegative whenever
 $R_{\mi}\le 0,$ which we have assumed from the beginning of the
 section. Let $\bar{R}$ denote the limit of $\hat{R}(t)$ as
 $t\to\infty.$
 \par Assume for a moment that $\bar{V}>0.$ Then it follows from
 (7.2) and (7.3) that $R_{\mi}(t)$ is asymptotic to
 $-\frac{3}{2t};$ in other words,
 $\bar{R}\bar{V}^{-\frac{2}{3}}=-\frac{3}{2}.$ Now the inequality
 (7.4) implies that whenever we have a sequence of parabolic
 neighborhoods $P(x^{\al},t^{\al},r\sqrt{t^{\al}},-r^2t^{\al}),$ for
 $t^{\al}\to\infty$ and some fixed small $r>0,$ such that the
 scalings of our solution with factor $t^{\al}$ smoothly converge
 to some limit solution, defined in an abstract parabolic neighborhood
 $P(\bar{x},1,r,-r^2),$ then the scalar curvature of this limit
 solution is independent of the space variables and equals
 $-\frac{3}{2t}$ at time $t\in [1-r^2,1];$ moreover, the strong
 maximum principle for (7.1) implies that the sectional curvature
 of the limit at time t is constant and equals $-\frac{1}{4t}.$
  This conclusion is also valid without the a priori
 assumption that $\bar{V}>0,$ since otherwise it is vacuous.
\par Clearly the inequalities and conclusions above hold for the
solutions to the Ricci flow with $\de(t)$-cutoff, defined in the
previous sections. From now on we assume that we are given such a
solution, so the estimates below may depend on it.
\par {\bf 7.2 Lemma.} {\it (a) Given $w>0,r>0,\xi>0$ one can find
$T=T(w,r,\xi)<\infty,$ such that if the ball
$B(x_0,t_0,r\sqrt{t_0})$ at some time $t_0\ge T$ has volume at
least $wr^3$ and sectional curvature at least $-r^{-2}t_0^{-1},$
then curvature at $x_0$ at time $t=t_0$ satisfies
\begin{equation}|2tR_{ij}+g_{ij}|< \xi. \end{equation}
\par (b) Given in addition $A<\infty$ and allowing $T$ to depend
on $A,$ we can ensure (7.5) for all points in
$B(x_0,t_0,Ar\sqrt{t_0}).$
\par (c) The same is true for $P(x_0,t_0,Ar\sqrt{t_0},Ar^2t_0).$}
\par {\it Proof.} (a) If $T$ is large enough then we can apply
corollary 6.8 to the ball $B(x_0,t_0,r_0)$ for
$r_0=\mi(r,\bar{r}(w))\sqrt{t_0};$ then use the conclusion of 7.1.
\par (b) The curvature control in $P(x_0,t_0,r_0/4,-\tau r_0^2),$
provided by corollary 6.8, allows us to apply proposition 6.3
(a),(b) to a controllably smaller neighborhood
$P(x_0,t_0,r_0',-(r_0')^2).$ Thus  by 6.3(b) we know that each
point in $B(x_0,t_0,Ar\sqrt{t_0})$ with scalar curvature at least
$Q=K_1'(A)r_0^{-2}$ has a canonical neighborhood. This implies
that for $T$ large enough such points do not exist, since if there
was a point with $R$ larger than $ Q,$ there would be a point
having a canonical neighborhood with $R=Q$ in the same ball, and
that contradicts the already proved assertion (a). Therefore we
have curvature control in the ball in question, and applying
6.3(a) we also get volume control there, so our assertion has been
reduced to (a).
\par (c) If $\xi$ is small enough, then the solution in the ball
$B(x_0,t_0,Ar\sqrt{t_0})$ would stay almost homothetic to itself
on the time interval $[t_0,t_0+Ar^2t_0]$ until (7.5) is violated
at some (first) time $t'$ in this interval. However, if  $T$ is
large enough, then this violation could not happen, because we can
apply the already proved assertion (b) at time $t'$ for somewhat
larger $A.$
\par {\bf 7.3} Let $\rho(x,t)$ denote the radius $\rho$ of the
ball $B(x,t,\rho)$ where $\mathrm{inf}\ Rm=-\rho^{-2}.$ It follows
from corollary 6.8, proposition 6.3(c), and the pinching estimate
(5.1) that for any $w>0$ we can find $\bar{\rho}=\bar{\rho}(w)>0,$
such that if $\rho(x,t)<\bar{\rho}\sqrt{t},$ then \begin{equation}
Vol\ B(x,t,\rho(x,t))<w\rho^3(x,t), \end{equation} provided that
$t$ is large enough (depending on $w$).
\par Let $M^-(w,t)$ denote the thin part of $M,$ that is the set
of $x\in M$ where (7.6) holds at time $t,$ and let $M^+(w,t)$ be
its complement. Then for $t$ large enough (depending on $w$) every
point of $M^+$ satisfies the assumptions of lemma 7.2. \par Assume
first that for some $w>0$ the set $M^+(w,t)$ is not empty for a
sequence of $t\to\infty.$ Then the arguments of Hamilton [H 4,\S
8-12] work in our situation. In particular, if we take a sequence
of points $x^{\al}\in M^+(w,t^{\al}), \ t^{\al}\to\infty,$ then
the scalings of $g_{ij}^{\al}$ about $x^{\al}$ with factors
$(t^{\al})^{-1}$ converge, along a subsequence of $\al\to\infty,$
to a complete hyperbolic manifold of finite volume. The limits may
be different for different choices of $(x^{\al},t^{\al}).$  If
none of the limits is closed, and $H_1$ is such a limit with the
least number of cusps, then, by an argument in [H 4,\S 8-10],
based on hyperbolic rigidity, for all sufficiently small $w',\
0<w'<\bar{w}(H_1),$ there exists a standard truncation $H_1(w')$
of $H_1,$ such that, for $t$ large enough, $M^+(w'/2,t)$ contains
an almost isometric copy of $H_1(w'),$ which in turn contains a
component of $M^+(w',t);$ moreover, this embedded copy of
$H_1(w')$ moves by isotopy as $t$ increases to infinity. If for
some $w>0$ the complement $M^+(w,t)\setminus H_1(w)$ is not empty
for a sequence of $t\to\infty,$ then we can repeat the argument
and get another complete hyperbolic manifold $H_2,$ etc., until we
find a finite collection of $H_j, 1\le j\le i,$ such that for each
sufficiently small $w>0$ the embeddings of $H_j(w)$ cover
$M^+(w,t)$  for all sufficiently large $t.$
\par Furthermore, the boundary tori of $H_j(w)$ are incompressible
in $M.$ This is proved [H 4,\S 11,12] by a minimal surface
argument, using a result of Meeks and Yau. This argument  does not
use the uniform bound on the normalized curvature, and goes
through even in the presence of surgeries, because the area of the
least area disk in question can only decrease when we make a
surgery.
\par {\bf 7.4} Let us redefine the thin part in case the thick one
isn't empty, $\ {M}^-(w,t)=M\setminus(H_1(w)\cup...\cup H_i(w)).$
Then, for sufficiently small $w>0$ and sufficiently large $t,$ $\
\ M^-(w,t)$ is diffeomorphic to a graph manifold, as implied by
the following general result on collapsing with local lower
curvature bound, applied to the metrics $t^{-1}g_{ij}(t).$
\par {\bf Theorem.} {\it Suppose $(M^{\al},g_{ij}^{\al})$ is a
sequence of compact oriented riemannian 3-manifolds, closed or
with convex boundary, and $w^{\al}\to 0.$ Assume that
\par (1) for each point $x\in M^{\al}$ there exists a radius
$\rho=\rho^{\al}(x), 0<\rho<1,$ not exceeding the diameter of the
manifold, such that the ball $B(x,\rho)$ in the metric
$g_{ij}^{\al}$ has volume at most $w^{\al}\rho^3$ and sectional
curvatures at least $-\rho^{-2};$
\par (2) each component of the boundary of $M^{\al}$ has diameter at most
$w^{\al},$ and has a (topologically trivial) collar of length one,
where the sectional curvatures are between $-1/4-\e$ and
$-1/4+\e;$
\par (3) For every $w'>0$ there exist $\bar{r}=\bar{r}(w')>0$ and
$K_m=K_m(w')<\infty ,\ m=0,1,2...,$ such that if $\al$ is large
enough, $0<r\le\bar{r},$ and the ball $B(x,r)$ in $g_{ij}^{\al}$
has volume at least $w'r^3$ and sectional curvatures at least
$-r^2, $ then the curvature and its $m$-th order covariant
derivatives at $x, \ m=1,2...,$ are bounded by $K_0r^{-2}$ and
$K_mr^{-m-2}$ respectively.
\par Then $M^{\al}$ for sufficiently large ${\al}$ are
diffeomorphic to  graph manifolds.}
\par Indeed, there is only one exceptional case, not covered by
the theorem above, namely, when $M=M^-(w,t),$ and $\rho(x,t),$ for
some $x\in M,$ is much larger than the diameter $d(t)$ of the
manifold, whereas the ratio $V(t)/d^3(t)$ is bounded away from
zero. In this case, since by the observation after formula (7.3)
the volume $V(t)$ can not grow faster than $\co \cdot
t^{\frac{3}{2}},$ the diameter does not grow faster than $\co
\cdot \sqrt{t},$ hence if we scale our metrics $g_{ij}(t)$ to keep
the diameter equal to one, the scaled metrics would satisfy the
assumption (3) of the theorem above and have the minimum of
sectional curvatures tending to zero. Thus we can take a limit and
get a smooth solution to the Ricci flow with nonnegative sectional
curvature, but not strictly positive scalar curvature. Therefore,
in this exceptional case $M$ is diffeomorphic to a flat manifold.
\par The proof of the theorem above will be given in a separate
paper; it has nothing to do with the Ricci flow; its main tool is
the critical point theory for distance functions and maps, see
[P,\S 2] and references therein. The assumption (3) is in fact
redundant; however, it allows to simplify the proof quite a bit,
by avoiding
 3-dimensional Aleksandrov spaces, and in particular, the
non-elementary Stability Theorem. \par Summarizing, we have shown
that for large $t$ every component of the solution is either
diffeomorphic to a graph manifold, or to a closed hyperbolic
manifold, or can be split by a finite collection of disjoint
incompressible tori into parts, each being diffeomorphic to either
a graph manifold or to a complete noncompact hyperbolic manifold
of finite volume. The topology of graph manifolds is well
understood [W]; in particular, every graph manifold can be
decomposed in a connected sum of irreducible graph manifolds, and
each irreducible one can in turn be split by a finite collection
of disjoint incompressible tori into Seifert fibered manifolds.
\section{On the first eigenvalue of the operator $-4\triangle +R$  }
\par {\bf 8.1} Recall from [I,\S 1,2] that Ricci flow is the gradient flow
for the first eigenvalue $\la$ of the operator $-4\triangle +R;$
moreover, $\frac{d}{dt}\la(t)\ge \frac{2}{3}\la^2(t)$ and
$\la(t)V^{\frac{2}{3}}(t)$ is non-decreasing whenever it is
nonpositive. We would like to extend these inequalities to the
case of Ricci flow with $\de(t)$-cutoff. Recall that we
immediately remove components with nonnegative scalar curvature.
\par {\bf Lemma.} {\it Given any positive continuous function $\xi(t)$
one can chose $\de(t)$ in such a way that for any solution to the
Ricci flow with $\de(t)$-cutoff, with normalized initial data, and
any surgery time $T_0,$ after which there is at least one
component, where the scalar curvature is not strictly positive, we
have an estimate $\la^+(T_0)-\la^-(T_0)\ge
\xi(T_0)(V^+(T_0)-V^-(T_0)),$ where $V^-,V^+$ and $\la^-,\la^+$
are the volumes and the first eigenvalues of $-4\triangle +R$
before and after the surgery respectively.}
\par {\it Proof.} Consider the minimizer $a$ for the functional

\begin{equation} \int (4|\nabla a|^2+Ra^2)\end{equation} under normalization
$\int a^2=1,$ for the metric after the surgery on a component
where scalar curvature is not strictly positive. Clearly is
satisfies the equation \begin{equation} 4\triangle a=Ra-\la^-a
\end{equation} Observe that since the metric contains an $\e$-neck
of radius about $r(T_0), $ we can estimate $\la^-(T_0)$ from above
by about $r(T_0)^{-2}$.
\par Let $M_{cap}$ denote the cap, added by the surgery. It is
attached to a long tube, consisting of $\e$-necks of various
radii. Let us restrict our attention to a maximal subtube, on
which the scalar curvature  at each point is at least
$2\la^-(T_0).$ Choose any $\e$-neck in this subtube, say, with
radius $r_0,$ and consider the distance function with range
$[0,2\e^{-1}r_0],$ whose level sets $M_z$ are almost round
two-spheres; let $M_z^+\supset M_{cap}$ be the part of $M,$
chopped off by $M_z.$ Then $$ \int_{M_z}
-4aa_z=\int_{M_z^+}(4|\nabla
a|^2+Ra^2-\la^-a^2)>r_0^{-2}/2\int_{M_z^+}a^2 $$ On the other
hand,
$$ |\int_{M_z}2aa_z-(\int_{M_z}a^2)_z|\le \co \cdot \int_{M_z}\e
r_0^{-1} a^2
$$ These two inequalities easily imply that $$\int_{M_0^+}a^2\ge
\mathrm{exp}(\e^{-1}/10)\int_{M_{\e^{-1}r_0}^+}a^2$$ Now the
chosen subtube contains at least about
$-\e^{-1}\mathrm{log}(\la^-(T_0)h^2(T_0))$ disjoint $\e$-necks,
where $h$ denotes the cutoff radius, as before. Since $h$ tends to
zero with $\de,$ whereas $r(T_0),$ that occurs in the bound for
$\la^-,$ is independent of $\de,$ we can ensure that the number of
necks is greater then $\mathrm{log}\ h,$ and therefore,
$\int_{M_{cap}}a^2<h^6,$ say. Then standard estimates for the
equation (8.2) show that $|\nabla a|^2$ and $Ra^2$ are bounded by
$\co \cdot h$ on $M_{cap},$ which makes it possible to extend $a$
to the metric before surgery in such a way that the functional
(8.1) is preserved up to $\co \cdot h^4.$ However, the loss of
volume in the surgery is at least $h^3,$ so it suffices to take
$\de$ so small that $h$ is much smaller than $\xi.$
\par {\bf 8.2} The arguments above lead to the following result {\it \par
 (a) If $(M,g_{ij})$ has $\la>0,$ then, for an appropriate choice of the cutoff
 parameter, the solution becomes extinct in finite time. Thus, if $M$ admits
 a metric with $\la>0$ then it is diffeomorphic to a connected sum of a finite collection
 of $\mathbb{S}^2\times\mathbb{S}^1$ and metric quotients of the round
 $\mathbb{S}^3.$ Conversely, every such connected sum admits a
 metric with $R>0,$ hence with $\la>0.$
   \par (b) Suppose $M$ does not admit any metric with
$\la>0,$ and let $\bar{\la}$ denote the supremum of $\la
V^{\frac{2}{3}}$ over all metrics on this manifold. Then
$\bar{\la}=0$ implies that  $M$  is a graph manifold. Conversely,
a graph manifold can not have $\bar{\la}<0.$
\par (c) Suppose $\bar{\la}<0$ and let $\bar{V}=
(-\frac{2}{3}\bar{\la})^{\frac{3}{2}}.$ Then $\bar{V}$ is the
minimum of $V,$ such that $M$ can be decomposed in connected sum
of a finite collection of $\mathbb{S}^2\times\mathbb{S}^1,$ metric
quotients of the round $\mathbb{S}^3,$ and some other components,
the union of which will be denoted by  $M',$ and there exists a
(possibly disconnected) complete hyperbolic manifold, with
sectional curvature $-1/4$ and volume $V,$ which can be embedded
in $M'$ in such a way that the complement (if not empty) is a
graph manifold. Moreover, if such a hyperbolic manifold has volume
$\bar{V},$ then its cusps (if any) are incompressible in $M'.$ }
\par For the proof  one needs in addition easily verifiable statements
that one can put metrics on connected sums preserving the lower
bound for scalar curvature [G-L], that one can put metrics on
graph manifolds with scalar curvature bounded below and volume
tending to zero [C-G], and that one can close a compressible cusp,
preserving the lower bound for scalar curvature and reducing the
volume, cf. [A,5.2]. Notice that using these results we can avoid
the hyperbolic rigidity and minimal surface arguments, quoted in
7.3, which, however, have the advantage of not requiring any a
priori topological information about the complement of the
hyperbolic piece.
\par The results above are exact analogs of the conjectures for the
Sigma constant, formulated by Anderson [A], at least in the
nonpositive case.
\section*{References}
 \ \ \  [I] G.Perelman The entropy formula for the Ricci flow
and its geometric applications. arXiv:math.DG/0211159 v1
\par [A] M.T.Anderson Scalar curvature and geometrization
conjecture for three-manifolds. Comparison Geometry (Berkeley,
1993-94), MSRI Publ. 30 (1997), 49-82.
\par [C-G] J.Cheeger, M.Gromov Collapsing Riemannian manifolds
while keeping their curvature bounded I. Jour. Diff. Geom. 23
(1986), 309-346.
\par [G-L] M.Gromov, H.B.Lawson Positive scalar curvature and the
Dirac operator on complete Riemannian manifolds. Publ. Math. IHES
58 (1983), 83-196.
\par [H 1] R.S.Hamilton Three-manifolds with positive Ricci
curvature. Jour. Diff. Geom. 17 (1982), 255-306.
\par [H 2] R.S.Hamilton Formation of singularities in the Ricci
flow. Surveys in Diff. Geom. 2 (1995), 7-136.
\par [H 3] R.S.Hamilton The Harnack estimate for the Ricci flow.
Jour. Diff. Geom. 37 (1993), 225-243.
\par [H 4] R.S.Hamilton Non-singular solutions of the Ricci flow
on three-manifolds. Commun. Anal. Geom. 7 (1999), 695-729.
\par [H 5] R.S.Hamilton Four-manifolds with positive isotropic
curvature. Commun. Anal. Geom. 5 (1997), 1-92.
\par G.Perelman Spaces with curvature bounded below. Proceedings
of ICM-1994, 517-525.
\par F.Waldhausen Eine Klasse von 3-dimensionalen
Mannigfaltigkeiten I,II. Invent. Math. 3 (1967), 308-333 and 4
(1967), 87-117.

\end{document}